\documentstyle[12pt]{article}
\begin{document}

\def\Diff{\hbox{\rm Diff}}

\title{Three Cocycles on $\Diff(S^1)$ Generalizing the 
 Schwarzian Derivative }
\author{S. Bouarroudj \and V.Yu. Ovsienko} 
\date{}
\maketitle 
\vskip 0,25 cm

{\abstract{
The first group of differentiable cohomology of $\Diff(S^1)$,
vanishing on the M\"obius subgroup 
$PSL(2,{\bf R})\subset\Diff(S^1)$,
with coefficients
in modules of linear differential
operators on $S^1$
is calculated.
We introduce three non-trivial 
$PSL(2,{\bf R})$-invariant 1-cocycles on $\Diff(S^1)$
generalizing the Schwarzian derivative.
}}

\section{Introduction}

{\bf 1.1}.
Consider the group $\Diff({\bf RP}^1)$
of diffeomorphisms of the circle ${\bf RP}^1\cong S^1$.
The well-known expression: 
$$
S(f)=\frac{f'''}{f'}-\frac{3}{2}\left(\frac{f''}{f'}\right)^2,
$$
where $f=f(x)\in \Diff(S^1)$,
$x$ is the affine parameter on ${\bf RP}^1$ and $f'=df(x)/dx$,
is called the {\it Schwarzian derivative}.

\vskip 0,3cm

The main properties of the Schwarzian derivative are
as follows:

\vskip 0,3cm
\noindent
(a) It satisfies the relation:
$$
S(f\circ g)=
S(f)\circ g\cdot(g')^2+S(g).
$$
That means, $S$ is a {\it 1-cocycle} on $\Diff({\bf RP}^1)$
with values in the space ${\cal F}_2$
of quadratic differentials (see \cite{ki});

\vskip 0,3cm
\noindent
(b) The kernel of the cocycle $S$ is the group of M\"obius
(linear-fractional or projective)
transformations $PSL(2,{\bf R})\subset\Diff({\bf RP}^1)$:
$S(f)\equiv0$ if and only if $f\in PSL(2,{\bf R})$.
Since the Schwarzian derivative
is a 1-cocycle, this means that it is {\it projectively invariant}.

\vskip 0,3cm
\noindent
{\bf 1.2}
The following two remarks have already been known to classics 
(see e.g. \cite{wi},\cite{ca},\cite{car}\footnote{We are grateful
to B. Kostant for this reference}):

\vskip 0,3cm
\noindent
(c) Given a Sturm-Liouville equation:
$
2\psi''+u(x)\psi=0,
$
where $u(x)\in C^{\infty}({\bf RP}^1)$,
let $\psi_1$ and $\psi_2$ be two independent solutions,
the potential $u(x)$ can be expressed as a function of the quotient:
$u=S(\psi_1/\psi_2)$.

\vskip 0,3cm
\noindent
(d) The same fact in other words.
Consider the space of Sturm-Liouville operators:
$$
A_u=2\frac{d^2}{dx^2}+u(x).
$$
The natural action of the group $\Diff({\bf RP}^1)$
on the space of Sturm-Liouville operators is:
\begin{equation}
\label{fu}
f^{-1}(A_u)=A_{u\circ f\cdot(f')^2}+S(f)
\end{equation}
To define this action, one considers the arguments of the
Sturm-Liouville operators as $-1/2$-densities and their
images as $3/2$-densities on ${\bf RP}^n$
(see Section 2 for the details).

\vskip 0,3cm
\noindent
{\bf 1.3}
We calculate the first group of differentiable cohomology of
$\Diff({\bf RP}^1)$ with
coefficients in the modules of
linear differential operators on tensor-densities,
vanishing on $PSL(2,{\bf R})$.

The main result of this paper is an explicit construction of
three families of non-trivial cocycles on $\Diff({\bf RP}^1)$
generalizing the Schwarzian derivative.
These cocycles 
are with values in the space of
linear differential operators on ${\bf RP}^1$.
They satisfy the main property of the Schwarzian derivative:
$PSL(2,{\bf R})$-invariance.

\section{$\Diff({\bf RP}^1)$-module structures
on the space of differential operators}

Space of linear differential operators on a manifold 
considered as a module over the group of diffeomorphisms,
is a classical subject. 
In the one-dimensional case, we refer 
\cite{wi} and \cite{ca}.
This subject is closely related to quantization 
(cf. \cite{ko}).
The modules of linear differential operators
on the space of tensor-densities on a smooth manifold
was studied in a series of recent papers
(see \cite{do},\cite{lmt},\cite{lo},\cite{go},\cite{ga},\cite{ma}). 

In this section we recall the definition of the
natural two-parameter family of modules
over the group of diffeomorphisms on the space of linear differential
operators.

\vskip 0,3cm
\noindent
{\bf 2.1}
Consider a one-parameter family of $\Diff({\bf RP}^1)$-actions 
on $C^{\infty}({\bf RP}^1)$:
$$
f^*_{\lambda}(\phi)=
\phi\circ f^{-1}\cdot({f^{-1}}^{\prime})^{\lambda}
$$
Denote ${\cal F}_{\lambda}$ the defined
$\Diff({\bf RP}^1)$-module structure
on $C^{\infty}({\bf RP}^1)$.
It is called the {\it module
tensor-densities of degree} $\lambda$ on ${\bf RP}^1$.
We will use the standard notation:
$
\phi=\phi(x)(dx)^{\lambda}.
$

For example, 
according to (\ref{fu}), the potential $u$ 
of a Sturm-Liouville operator should be considered as an
element of ${\cal F}_2$ 
(so-called {\it quadratic differential}: $u=u(x)(dx)^2$,
see e.g. \cite{ki}).

\vskip 0,3cm
\noindent
{\bf 2.2}
Denote ${\cal D}^{k}$ the space of  $k^{th}$-order
differential operators 
\begin{equation}
\label{Op}
A=a_k(x)\frac{d^k}{dx^k}+\cdots +a_0(x)
\end{equation}
where $a_i(x)\in C^{\infty }({\bf RP}^1).$

\vskip 0,3cm
\noindent
{\bf Definition}.
A two-parameter family of actions of $\Diff({\bf RP}^1)$ on the
space of differential operators (\ref{Op}) is defined by:
\begin{equation}
\label{act}
g_{\lambda,\mu }(A)=g^*_{\mu}\circ A\circ(g^*_{\lambda})^{-1}
\end{equation}
In other words, we consider differential operators
acting on tensor-densities:
$A:{\cal F}_{\lambda}\to{\cal F}_{\mu}$.

\vskip 0,3cm

Denote by ${\cal D}^{k}_{\lambda,\mu}$ the space of operators
(\ref{Op}) endowed with the defined $\Diff({\bf RP}^1)$-module
structure.

\vskip 0,3cm
\noindent
{\bf Remark}. 
The complete
classification of modules ${\cal D}^{k}_{\lambda,\mu}$ was obtained in
\cite{ga}.

\vskip 0,3cm

In this paper we consider the action of the
group $\Diff({\bf RP}^1)$ on this space
and study the cohomology groups
arising in this context.

\section{The main result}

Studying of the modules of differential operators 
leads to the cohomology group
$H^1({\rm Diff}({\bf RP}^1);
{\rm Hom}({\cal F}_{\lambda},{\cal F}_{\mu}))$
 (see \cite{lo}).
This cohomology can be considered as a
measure of ``non-triviality'' for the modules of 
linear differential operators.

\vskip 0,3cm
\noindent
{\bf 3.1}. 
Let us formulate the precise problem considered in this paper.

We will study the {\it differentiable} (or local)
cohomology (see \cite{fu} for the general definition). 
This means, we consider
only differentiable cochains on $\Diff({\bf RP}^1)$ 
with coefficients in the space of differential operators:
${\cal D}_{\lambda,\mu}
\subset
{\rm Hom}({\cal F}_{\lambda},{\cal F}_{\mu})$.

We will impose one more condition:
$PSL(2,{\bf R})$-{\it invariance}. 
In other words,
we consider only the cohomology classes {\it vanishing} on
$PSL(2,{\bf R})$
(note that for the case of first cohomology these two notions coincide).

\proclaim Theorem 1.
The differentiable
cohomology of $\Diff({\bf RP}^1)$ with coefficients 
in the module of linear differential operators,
vanishing on $PSL(2,{\bf R})$:
$$
H^1_{\rm diff}({\rm Diff}({\bf RP}^1),PSL(2,{\bf R});\;
{\cal D}_{\lambda,\mu})
$$
is one-dimensional in the following cases:
\hfill\break
(a)
$\mu-\lambda=2, \lambda\not=-1/2$,
\hfill\break
(b)
$\mu-\lambda=3, \lambda\not=-1$, 
\hfill\break
(c)
$\mu-\lambda=4, \lambda\not=-3/2$, 
\hfill\break
(d)
$(\lambda,\mu)=(-4,1),(0,5)$.
\hfill\break
Otherwise, this cohomology group is trivial.
\par

We will prove this theorem in Section 6.

\vskip 0,3cm

Therefore, there exist three families of non-trivial
cohomology classes depending on $\lambda$ as a parameter
(and with fixed $\mu=\lambda+2,\lambda+3$ or $\lambda+4$).
We will call these cohomology classes {\it stable}.

\proclaim Corollary.
For the {\it generic} values of $\lambda$ one has:
$$
H^1_{\rm diff}({\rm Diff}({\bf RP}^1),PSL(2,{\bf R});\;
{\cal D}_{\lambda,\mu})=
\left\{
\matrix{
{\bf R}, \;\mu-\lambda=2,3,4\hfill\cr
0, \;\;{\rm otherwise}\hfill\cr
}
\right.
$$
\par

Let us give here the list of 1-cocycles generating the 
stable non-trivial
cohomology classes.

\vskip 0,3cm
\noindent
{\bf 3.2}. 
The main result of this paper is:

\proclaim Theorem 2. 
(i)
For every $\lambda$ there exist unique (up to a multiple)
1-cocycles:
$$
\matrix{
{\cal S}_{\lambda}:\Diff({\bf RP}^1)\to
{\cal D}_{\lambda,\lambda+2},\hfill\cr\noalign{\smallskip}
{\cal T}_{\lambda}:\Diff({\bf RP}^1)\to
{\cal D}_{\lambda,\lambda+3},\hfill\cr\noalign{\smallskip}
{\cal U}_{\lambda}:\Diff({\bf RP}^1)\to
{\cal D}_{\lambda,\lambda+4}
}
$$
vanishing on $PSL(2,{\bf R})$. They are
given by the formul{\ae}:
\begin{equation}
\label{sch}
\matrix{
{\cal S}_{\lambda}(f)\;\;\;=\;\;\;S(f)
\hfill\cr\noalign{\bigskip}
\displaystyle{\cal T}_{\lambda}(f)\;\;\;=\;\;\;S(f)\frac{d}{dx}-
\frac{\lambda}{2} S(f)^{\prime}
\hfill\cr\noalign{\bigskip}
\matrix{{\cal U}_{\lambda}(f)&=&\displaystyle S(f)\frac{d^2}{dx^2}-
\frac{2\lambda+1}{2} S(f)^{\prime}\frac{d}{dx}
\hfill\cr\noalign{\smallskip}
&&\displaystyle+\frac{\lambda(2\lambda+1)}{10} S(f)^{\prime\prime}-
\frac{\lambda(\lambda+3)}{5} S(f)^2,
}
}
\end{equation}
where $S(f)$ is the Schwarzian derivative.
\hfill\break
(ii)
The cocycles ${\cal S}_{\lambda},{\cal T}_{\lambda}$ and ${\cal U}_{\lambda}$
are non-trivial for every $\lambda$ except
$\lambda=-1/2,\lambda=-1$ and $\lambda=-3/2$ respectively.
\par

Explicit formul{\ae} for the cocycles with values in the exceptional
modules ${\cal D}_{-4,1}$ and ${\cal D}_{0,5}$ (cf. Theorem 1)
will be given in the end of this paper.

\vfill\eject
\noindent
{\bf 3.3. Proof of existence}. Let us  show that the maps
${\cal S}_{\lambda},{\cal T}_{\lambda}$ and ${\cal U}_{\lambda}$
are 1-cocycles.

(1)
The first map ${\cal S}_{\lambda}(f)$ is just a zero-order 
differential operator of multiplication by the Schwarzian derivative:
${\cal S}_{\lambda}(f)(\phi)=S(f)\cdot\phi$.
The condition of 1-cocycle for ${\cal S}_{\lambda}$
follows from those for $S$.

(2)
In general, if $k\geq2$, the modules ${\cal D}^{k}_{\lambda,\mu}$
are not isomorphic to the modules of tensor-densities.
However, for $k=1$
this is steel the case. The result is as follows:

$$
 {\rm if}\quad \mu-\lambda\not=1,\quad {\rm then}\quad
{\cal D}^1_{\lambda,\mu}\cong
{\cal F}_{\mu-\lambda-1}\oplus{\cal F}_{\mu-\lambda}.
$$
Given an operator 
$A=a_1(x)\frac{d}{dx}+a_0(x)\in{\cal D}^1_{\lambda,\mu}$,
the principal symbol $a_1$ transforms under the action
(\ref{act}) as a $\mu-\lambda-1$-density.
Verify that the quantity
$\bar a_0:=a_0-\frac{\lambda}{\mu-\lambda-1}a_1'$
transforms as a $\mu-\lambda$-density.
The isomorphism is as follows:
$$
\sigma(A)=
\left(a_1(x)(dx)^{\mu-\lambda-1},
\bar a_0(x)(dx)^{\mu-\lambda}\right).
\eqno{(\ref{sym}')}
$$

Existence of the cocycle ${\cal T}_{\lambda}$
is a corollary of the isomorphism $\sigma$.
Namely, the map 
$\overline{\cal T}_{\lambda}=\sigma\circ{\cal T}_{\lambda}:
\Diff({\bf RP}^1)\to
{\cal F}_2\oplus{\cal F}_3$
is given by the formula:
$$
\overline{\cal T}_{\lambda}(f)=\left(S(f)(dx)^2,0\right).
\eqno{(\ref{sch}')}
$$
This is obviously a 1-cocycle.

(3) Let us show that the map ${\cal U}_{\lambda}$ is
a 1-cocycles for every $\lambda$.

The modules of second order differential operators
${\cal D}^2_{\lambda,\lambda+3}$ is not isomorphic
to a direct sum of tensor-densities, for every
$\lambda$ except
$\lambda=0,-3$.
However, if $\mu-\lambda\not=1,3/2,2$ there exists a linear map
$$
\sigma:{\cal D}^2_{\lambda,\mu}\to
{\cal F}_{\mu-\lambda-2}\oplus
{\cal F}_{\mu-\lambda-1}\oplus
{\cal F}_{\mu-\lambda},
$$
which is $PSL(2,{\bf R})$-equivariant
(cf.\cite{cmz},\cite{lo},\cite{ga} and Section 7).
This defines a symbol of a differential operator
$A=a_2\frac{d^2}{dx^2}+a_1\frac{d}{dx}+a_0\in{\cal D}^2_{\lambda,\mu}$
in a canonical way:
$$
\sigma(A)=
\left(a_2(dx)^{\mu-\lambda-2},\;
\bar a_1(dx)^{\mu-\lambda-1},\;
\bar a_0(dx)^{\mu-\lambda}\right),
\eqno{(\ref{sym}'')}
$$
where 
$$
\bar a_1=a_1+\frac{2\lambda+1}{\mu-\lambda-2}a_2'\quad
{\rm and}\quad
\bar a_0=a_0+
\frac{\lambda}{\mu-\lambda-1}a_1+
\frac{\lambda(2\lambda+1)}{(\mu-\lambda)(2(\mu-\lambda)-3)}a_2''.
$$
These expressions is a partial case of the formula
(\ref{sym}) below.

One can easily calculate
the $\Diff({\bf RP}^1)$-action (\ref{act}) 
on the right hand side:
$$
\sigma\circ f_{\lambda,\mu}\circ\sigma^{-1}:
(a_2,\bar a_1,\bar a_0)\mapsto
\left(f^*(a_2),
f^*(a_1),
f^*(a_0)+\beta S(f^{-1})f^*(a_2)\right),
\eqno{(\ref{act}')}
$$
where $\beta=\frac{2\lambda(\mu-1)}{2(\mu-\lambda)-3}$.
In the particular case: $\mu-\lambda=4$
one obtains: 
$$
\beta=\frac{2\lambda(\lambda+3)}{5}.
$$

Verify that the map
$\overline{\cal U}_{\lambda}=
\sigma\circ{\cal U}_{\lambda}$ is of the form:
$$
\overline{\cal U}_{\lambda}(f)=
\left(S(f)(dx)^2,\;
0,\;
-\frac{\lambda(\lambda+3)}{5}S(f)^2(dx)^4\right).
\eqno{(\ref{sch}'')}
$$
Now, it is very easy to check that 
this map is a 1-cocycle with respect to the action
(\ref{act}').

It is proven, that the maps 
${\cal S}_{\lambda},{\cal T}_{\lambda}$ and ${\cal U}_{\lambda}$
are indeed 1-cocycles.

\vskip 0,3cm

Part (ii) of Theorem 2 and the
uniqueness of the cocycles (\ref{sch}) will be proven in 
the next section.

\vskip 0,3cm
\noindent
{\bf Remarks}. 1) Recall, that there exist two more analogues
of the Schwarzian derivative (see \cite{fu}):
$$
f\mapsto f'(x)
\quad {\rm and} \quad
f\mapsto \frac{f''(x)}{f'(x)}dx
$$
with values in ${\cal F}_0$ and ${\cal F}_1$ respectively
These 1-cocycles, however, do not vanish on $PSL(2,{\bf R})$.
The 1-cocycles with values in
${\cal D}_{\lambda,\lambda}$
and
${\cal D}_{\lambda,\lambda+1}$
defined through multiplication by these 1-cocycles
turn to be trivial for every
$\lambda\not=0$.

\vskip 0,3cm

2)
The module of second order differential operators
${\cal D}^2_{\lambda,\mu}$ is a direct sum of modules of tensor-densities:
$
{\cal D}^2_{\lambda,\mu}\cong
{\cal F}_{\mu-\lambda-2}\oplus
{\cal F}_{\mu-\lambda-1}\oplus
{\cal F}_{\mu-\lambda},
$
if and only if $\lambda=0$ or $\mu=1$
(cf. Section 7).
If $\lambda-\mu=3$, these conditions
corresponds to the case when the 1-cocycle
${\cal U}_{\lambda}$ is a {\it linear} map 
(namely, $\lambda=0$ or $\lambda=-3$).

\vskip 0,3cm

3) The unique module of third-order operators 
isomorphic to a direct sum of modules of
tensor-densities is
$
{\cal D}^3_{0,1}\cong
{\cal F}_{-2}\oplus
{\cal F}_{-1}\oplus
{\cal F}_0\oplus
{\cal F}_1.
$

\section{$PSL(2,{\bf R})$-invariant differential operators \allowbreak
and
trivialization of the cocycles 
${\cal S}_{\lambda},{\cal T}_{\lambda}$ \break
and ${\cal U}_{\lambda}$}

{\bf 4.1. Proof of Theorem 2, Part (ii)}.
Let ${\cal C}:\Diff({\bf RP}^1)\to{\cal D}_{\lambda,\mu}$
be a 1-cocycle vanishing on $PSL(2,{\bf R})$.
Suppose that ${\cal C}$ is a coboundary:
${\cal C}=\delta(B)$ for some $B\in{\cal D}_{\lambda,\mu}$.
This means, for every $f\in\Diff({\bf RP}^1)$,
$$
{\cal C}(f)=f_{\lambda,\mu}(B)-B.
$$
In particular, for $f\in PSL(2,{\bf R})$,
it follows that the {\it operator $B$ is
$PSL(2,{\bf R})$-equivariant}.

\vskip 0,3cm

$PSL(2,{\bf R})$-equivariant linear differential
operators on tensor-densities were classified in
\cite{bo}. Let us
recall here the classical result:

\proclaim Bol's Theorem. (i) For every $k$ there exists a unique (up to a
constant)
$PSL(2,{\bf R})$-invariant linear differential
operator of order $k$:
$
\partial^k:
{\cal F}_{\frac{1-k}{2}}\to
{\cal F}_{\frac{1+k}{2}},
$
such that
$$
\partial^k\left(\phi(x)(dx)^{\frac{1-k}{2}}\right)=
\phi^{(k)}(x)(dx)^{\frac{1+k}{2}}.
$$
\hfill\break
(ii) If $(\lambda,\mu)\not=((1-k)/2,(1+k)/2)$,
then there is no non-zero $PSL(2,{\bf R})$-invariant
linear differential operator from 
${\cal F}_{\lambda}$ to ${\cal F}_{\mu}$.
\par

The operator $\partial^k$ is called the {\it Bol operator}.

\vskip 0,3cm

Let us show how Theorem 2, Part (ii) follows from the Bol Theorem.

\vskip 0,3cm

(1) Verify that the cocycles 
${\cal S}_{-\frac{1}{2}},{\cal T}_{-1}$ and ${\cal U}_{-\frac{3}{2}}$
are differentials of the Bol operators:
$$
{\cal S}_{-\frac{1}{2}}=
\delta\left(\partial^2\right),\;\;
{\cal T}_{-1}=
\delta\left(\partial^3\right),\;\;
{\cal U}_{-\frac{3}{2}}=
\delta\left(\partial^4\right).
$$

\vskip 0,3cm

(2) On the other hand, Bol's Theorem implies that the cocycles 
${\cal S}_{\lambda},{\cal T}_{\lambda}$ and ${\cal U}_{\lambda}$
with $\lambda\not=-1/2,-1$ and $-3/2$ respectively,
are non-trivial.
Indeed, these cocycles vanish on $PSL(2,{\bf R})$.
Therefore, if these cocycles are trivial, they have to be
differential of some $PSL(2,{\bf R})$-invariant operators.
This contradicts to the Bol classification.

\vskip 0,3cm

Theorem 2, Part (ii) is proven.

\vskip 0,3cm
\noindent
{\bf 4.2. Proof of uniqueness}.
Let us show that uniqueness of the cocycles
${\cal S}_{\lambda},{\cal T}_{\lambda}$ and ${\cal U}_{\lambda}$
follows from Theorem 1 and the Bol Theorem.

Let $C_1$ and $C_2$ be two differentiable 1-cocycles on
$\Diff({\bf RP}^1)$
vanishing on $PSL(2,{\bf R})$, with values in the same module
${\cal D}_{\lambda,\mu}$.
We will show that $C_1$ and $C_2$ are proportional to each other
(doesn't matter whether $C_1$ and $C_2$ are trivial or not).

Since the considered cohomology group is at most one-dimensional,
there exists a linear combination 
$C=\alpha_1C_1+\alpha_2C_2$ which is a coboundary.
Since the cocycle $C$ vanishes on $PSL(2,{\bf R})$,
this means, $C=\delta(B)$, where $B\in{\cal D}_{\lambda,\mu}$
is a $PSL(2,{\bf R})$-invariant operator.
In the case $\lambda\not=-1/2,-1,-3/2,\dots$,
Bol's Theorem implies: $C\equiv0$.
If $\lambda=-1/2,-1,-3/2,\dots$, then the cohomology
is trivial. In this case Bol's Theorem implies that
both of the cocycles $C_1$ and $C_2$ are proportional
to the differential of one of the Bol operators.

\vskip 0,3cm

Theorem 2 is proven.

\section{Exceptional modules of differential ope\-rators}

The exceptional values
$\lambda=-1/2,-1,-3/2$
correspond to the particular modules 
of differential operators studied by classics.
In this section we interpret these modules
in terms of the constructed cocycles.

\vskip 0,3cm

Consider the space of differential operators
of the form:
\begin{equation}
\label{cla}
A=
\frac{d^k}{dx^k}+
a_{k-2}\frac{d^{k-2}}{dx^{k-2}}+
a_{k-3}\frac{d^{k-3}}{dx^{k-3}}+
\cdots+a_0,
\end{equation}
where $A$ is acting on the space of $(1-k)/2$-densities
with values in $(1+k)/2$-densities:
$$
A\in{\cal D}_{\frac{1-k}{2},\frac{1+k}{2}}.
$$
The structure of $\Diff({\bf RP}^1)$-module on
this space has already been
studied in \cite{wi} and \cite{ca}.

Note that the modules (\ref{cla}) are closely related to
so-called Adler-Gelfand-Dickey (or classical $W$)
Poisson structure.

\vskip 0,3cm
\noindent
{\bf Example (a)}. The module of Sturm-Liouville operators 
considered in the introduction
is a submodule ${\cal D}^2_{-\frac{1}{2},\frac{3}{2}}$.
Indeed, verify that
the formula of $\Diff({\bf RP}^1)$-action (\ref{fu})
coincides with the action $f_{-\frac{1}{2},\frac{3}{2}}$.
In the case $u\equiv0$,
one has:
$$
f_{-\frac{1}{2},\frac{3}{2}}\left(\partial^2\right)-
\partial^2=
{\cal S}_{-\frac{1}{2}}(f).
$$
This is just the coboundary relation:
${\cal S}_{-\frac{1}{2}}=\delta(\partial^2)$.

\vskip 0,3cm
\noindent
{\bf Example (b)}. The operators (\ref{cla})
for $k=3$ can be written in the form:
$$
A_{u,v}=
\frac{d^3}{dx^3}
+4u(x)\frac{d}{dx}
+2u(x)^{\prime}+v(x),
$$
where $A\in{\cal D}_{\lambda,\lambda+3}$.
This
is a particular case of this module corresponding to $k=2$.
It is easy to check that
the formula of $\Diff({\bf RP}^1)$-action reads as:
$$
g^{-1}_{\lambda,\lambda+3}(A_{u,v})=
A_{u\circ g\cdot(g')^2,v\circ g\cdot(g')^3}
+{\cal T}_{-1}(g).
$$
In the case $u=v\equiv0$ one obtains the relation
${\cal T}_{-1}=\delta(\partial^3)$.

\vskip 0,3cm
\noindent
{\bf Example (c)}. Every operator (\ref{cla})
for $k=4$ can be written in the form:
$$
A_{u,v,w}=
\frac{d^3}{dx^3}
+5u(x)\frac{d^2}{dx^2}
+5u(x)^{\prime}\frac{d}{dx}
+\frac{3}{2}u(x)^{\prime\prime}+
\frac{9}{4}u(x)^2
+v(x)\frac{d}{dx}
+w(x),
$$
The $\Diff({\bf RP}^1)$-action reads as:
$$
g^{-1}_{\lambda,\lambda+4}(A_{u,v,w})=
A_{u\circ g\cdot(g')^2,v\circ g\cdot(g')^3,w\circ g\cdot(g')^4}
+{\cal U}_{-\frac{3}{2}}(g).
$$
In the case $u=v\equiv0$ one obtains the relation
${\cal U}_{-\frac{3}{2}}=\delta(\partial^4)$.

Note, that the coefficients of the terms
containing $u(x)$ coincide with the coefficients
of the cocycle ${\cal U}_{\lambda}$ for 
$\lambda=-3/2$.

\vskip 0,3cm
\noindent
{\bf Remark}.
In the same way one can interpret the modules of operators
(\ref{cla}) for an arbitrary value of $k$, in terms of
1-cocycles defined as differentials of
the Bol operators.

\section{
Bilinear differential $PSL(2,{\bf R})$-invariant
\allowbreak
operators and 
cohomology of vector fields Lie algebra
}

In this section we prove Theorem 1.

Consider the differentiable cohomology 
of the Lie algebra of vector fields
vanishing on the M\"obius subalgebra 
$sl(2,{\bf R})\subset{\rm Vect}({\bf RP}^1)$:
\begin{equation}
\label{hom}
H^1_{\rm diff}({\rm Vect}({\bf RP}^1),sl(2,{\bf R});
{\cal D}_{\lambda,\mu}).
\end{equation}
Every class of differentiable cohomology
$H^1_{\rm diff}(\Diff({\bf RP}^1),PSL(2,{\bf R});
{\cal D}_{\lambda,\mu})$ 
corresponds to
some non-trivial class of ${\rm Vect}({\bf RP}^1)$-cohomology (\ref{hom})
(cf. \cite{fu}).

To prove Theorem 1, let us first show that
the similar result holds for the cohomology group (\ref{hom}).

\proclaim Proposition 6.1. The differentiable
cohomology (\ref{hom})
is one-dimensional in the following cases:
\hfill\break
(a)
$\mu-\lambda=2, \lambda\not=-1/2$,
\hfill\break
(b)
$\mu-\lambda=3, \lambda\not=-1$, 
\hfill\break
(c)
$\mu-\lambda=4, \lambda\not=-3/2$, 
\hfill\break
(d)
$(\lambda,\mu)=(-4,1),(0,5)$.
\hfill\break
Otherwise, this cohomology group is trivial.
\par

\vskip 0,3cm
\noindent
{\bf 6.1}. 
The first remark is as follows.

\proclaim Lemma 6.1.
Given a differentiable 1-cocycle $c$ on ${\rm Vect}({\bf RP}^1)$
vanishing on $sl(2,{\bf R})$,
with values in ${\cal D}_{\lambda,\mu}$,
the bilinear differential operator
$J:{\rm Vect}({\bf RP}^1)\otimes
{\cal F}_{\lambda}\to
{\cal F}_{\mu}$
defined by:
$$
J(X,\phi)=
c(X)(\phi),
$$
is $sl(2,{\bf R})$-invariant.
\par

\noindent
{\bf Proof}. 
Since $c$ a 1-cocycle, it satisfies the relation:
$$
L_X\circ c(Y)-c(Y)\circ L_X
-
L_Y\circ c(X)+c(X)\circ L_Y
=
c([X,Y])
$$
for every $X,Y\in{\rm Vect}({\bf RP}^1)$,
then for $X\in sl(2,{\bf R})$ one gets:
$$
L_X(J(Y,\phi))=
J([X,Y],\phi)+J(Y,L_X(\phi)).
$$
That means, $J$ is $sl(2,{\bf R})$-invariant.

\vskip 0,3cm
\noindent
{\bf 6.2}.
All the $sl(2,{\bf R})$-invariant bilinear differential
operators on tensor-densities:
${\cal F}_{\lambda}\otimes{\cal F}_{\mu}\to{\cal F}_{\nu}$
were classified in \cite{go} 
(see also \cite{gara} for a clear and detailed exposition).
The result is as follows. 

\proclaim Gordan's Theorem.
(i) For every $\lambda,\mu$ and integer $m\geq0$, there
exists a $sl(2,{\bf R})$-invariant bilinear differential operator
$J_m^{\lambda,\mu}:
{\cal F}_{\lambda}\otimes
{\cal F}_{\mu}\to
{\cal F}_{\lambda+\mu+m}$ given by:
$$
J_m^{\lambda,\mu}(\phi, \psi ) 
= 
\sum_{i+j=m} (-1)^i m!
{2\lambda+m-1 \choose i} {2\mu+m-1 \choose j}  
{\phi}^{(i)} {\psi}^{(j)} 
$$
\hfill\break
(ii) If either 
$\lambda$ or $\mu\not\in\{-1/2,-1,-3/2,\dots\}$, then
the operator $J_m^{\lambda,\mu}$ is the
unique (up to a constant)
$sl(2,{\bf R})$-invariant bilinear differential
operator from 
${\cal F}_{\lambda}\otimes{\cal F}_{\mu}$
to
${\cal F}_{\lambda+\mu+m}$.
\par

The operator $J_m^{\lambda,\mu}$ is called the {\it transvectant}.

\vskip 0,3cm

Therefore, one obtain a series of bilinear $sl(2,{\bf R})$-invariant
maps: 
$$
J_m^{-1,\lambda}:{\rm Vect}({\bf RP}^1)\otimes
{\cal F}_{\lambda}\to
{\cal F}_{\lambda+m-1}.
$$
If $\lambda\not=-1/2,-1,-3/2,\dots$,
the operator $J_m^{-1,\lambda}$ is a
{\it unique} (up to a constant) operator
invariant with respect to $sl(2,{\bf R})$.

\vskip 0,3cm

An important property of the operators $J_m^{-1,\lambda}$
with $m\geq3$ is that
they vanish on the subalgebra
$sl(2,{\bf R})$.

\vskip 0,3cm
\noindent
{\bf 6.3 Proof of Proposition 6.1}.

(1) Let us consider the case: 
$\lambda\not==-1/2,-1,-3/2,\dots$.

Every differentiable 1-cocycle $c$ on 
${\rm Vect}({\bf RP}^1)$
vanishing on $sl(2,{\bf R})$,
with values in ${\cal D}_{\lambda,\mu}$,
is proportional to the map $c_m$ with $m\geq3$ such that:
$$
c_m(X)(\phi):=J_m^{-1,\lambda}(X,\phi).
$$
Indeed, it follows from Lemma 6.2 that $c$
is $sl(2,{\bf R})$-invariant and it follows from 
the Gordan Theorem that it is proportional to
$J_m^{-1,\lambda}$.
Now,
let us check for each of the maps
$c_m$ if it is a 1-cocycles.

\vskip 0,3cm

It is easy to calculate their explicit formul{\ae}:
$$
\matrix{
\displaystyle J_3^{-1,\lambda}(X\frac{d}{dx},\phi)=
X'''\phi,\hfill\cr
\noalign{\smallskip}
\displaystyle J_4^{-1,\lambda}(X\frac{d}{dx},\phi)=
X'''\phi'-\frac{\lambda}{2}X^{IV}\phi\hfill\cr 
\noalign{\smallskip}
\displaystyle J_5^{-1,\lambda}(X\frac{d}{dx},\phi)=
X'''\phi''-\frac{2\lambda+1}{2}X^{IV}\phi'+
\frac{\lambda(2\lambda+1)}{10}X^{V}\phi\hfill\cr 
}
$$
These maps are the 1-cocycles on
${\rm Vect}({\bf RP}^1)$
corresponding to the cocycles
${\cal S}_{\lambda},{\cal T}_{\lambda}$ and ${\cal U}_{\lambda}$.
They are non-trivial 
if and only if $\lambda\not=-1/2,-1,-3/2$ respectively.

One can verify by direct
calculation that:

\proclaim Lemma 6.2.
The map $J_6^{-1,\lambda}$
defines a 1-cocycle if and only if $\lambda=-4,0,-2$.
\par

The first two cocycles define non-zero cohomology
classes of (\ref{hom}),
in the last case this cocycle is trivial 
(equals to the differential of the Bol operator).

\proclaim Lemma 6.3.
The map $J_m^{-1,\lambda}$ with $m\geq7$
defines a 1-cocycle if and only if $\lambda=(1-m)/2$
\par

The corresponding 1-cocycle is trivial.

\vskip 0,3cm

Proposition 6.1 is proven for the case 
$\lambda\not=-1/2,-1,-3/2,\dots$.

\vskip 0,3cm

(2) Let $\lambda\in\{-1/2,-1,-3/2,\dots\}$.
In this case,
the property of $sl(2,{\bf R})$-invariance
does not define a unique operator from
${\rm Vect}({\bf RP}^1)\otimes
{\cal F}_{\lambda}$ to
${\cal F}_{\lambda+m-1}.$
However, the two properties:
of $sl(2,{\bf R})$-invariance
and vanishing on $sl(2,{\bf R})$ determine the transvectants
$J_m^{-1,\lambda}$ uniquely.

\proclaim Lemma 6.4.
The operators $J_m^{-1,\lambda}$ with $m\geq2$ are the
unique (up to a constant) $sl(2,{\bf R})$-invariant bilinear differential
operators vanishing on $sl(2,{\bf R})$.
\par

\vskip 0,3cm
\noindent
{\bf Proof}: straightforward  (see \cite{lo}).

\vskip 0,3cm

Now, Proposition 6.1 follows from
Lemmas 6.2 and 6.3
are trivial for $\lambda=-1/2,-1,-3/2$ respectively.

\vskip 0,3cm

This gives an upper boundary for the dimension
of the cohomology group. 
After that, Theorem 1 follows from the explicit
construction of the cocycles generating non-trivial cohomology
classes.

\vskip 0,3cm

Theorem 1 is proven.

\vskip 0,3cm
\noindent
{\bf 6.4 Remark}.
Another way to prove Theorem 1 is to use the results
of Feigin-Fuchs \cite{ff} and Roger \cite{ro}.
The cohomology group 
$H^1(W;{\rm Hom}(F_{\lambda},F_{\mu}))$,
where $W$ is the Lie algebra of formal vector fields on $\bf R$
and $F_{\lambda}$ is a module of formal
tensor-densities on $\bf R$,
was calculated in \cite{ff}.
The analogous results was obtained in \cite{ro} 
in the differentiable case.
One can obtain Theorem 1 from the result of \cite{ff}
and \cite{ro} by selecting the cohomology 
classes trivial on $sl(2,{\bf R})$.

\section{Relations to the modules of 
higher order differential operators}

What follows is an illustration of what has been done.
We will show how are the cocycles
${\cal S}_{\lambda},{\cal T}_{\lambda}$ and ${\cal U}_{\lambda}$
appear in the modules of higher order differential operators.
In this section we will not give the details of calculations.

\vskip 0,3cm
\noindent
{\bf 7.1}.
Let us recall the following result from
 \cite{cmz} (see also \cite{lo} and \cite{ga})
concerning the restriction of the ${\rm Diff}({\bf RP}^1)$-module
structure to the subgroup $PSL(2,{\bf R})$:

\proclaim Cohen-Manin-Zagier's Theorem. For
$\rho-\nu\not=1,\frac{3}{2},2,\dots,k$, there exists an isomorphism of 
$PSL(2,{\bf R})$-modules:
$$
\sigma:{\cal D}^k_{\nu,\rho}\to
{\cal F}_{\rho-\nu-k}\oplus
{\cal F}_{\rho-\nu-k+1}
\oplus\cdots\oplus
{\cal F}_{\rho-\nu}
$$
\par

This isomorphism is
called the $PSL(2,{\bf R})$-{\it equivariant symbol map.}

\vskip 0,3cm

The explicit formula of the $PSL(2,{\bf R})$-equivariant symbol map is:
\begin{equation}
\label{sym}
\bar a_i(x)=
\sum_{j=i}^{k}\alpha_i^ja^{(j-i)}_j(x)
\quad\in{\cal F}_{\rho-\nu+i},
\end{equation}
where the constants $\alpha_i^j$ are written in terms of 
binomial coefficients:
$$
\alpha_i^j=
\frac{{j \choose i}{2\nu+i\choose
2\nu+j}}{{2(\rho-\nu)-i-j-1\choose 2(\rho-\nu)-2i-1}}
$$

Note, that the isomorphisms (\ref{sym}') and (\ref{sym}'')
are particular cases of the $PSL(2,{\bf R})$-equivariant symbol map.

\vskip 0,3cm
\noindent
{\bf 7.2}.
Now, consider the action (\ref{act}) of the group
${\rm Diff}({\bf RP}^1)$ on the modules ${\cal D}^k_{\lambda,\mu}$
with $\rho-\nu\not=1,\frac{3}{2},2,\dots,k$,
written in terms of $PSL(2,{\bf R})$-equivariant symbol:
$\sigma\circ f_{\nu,\rho}\circ\sigma^{-1}$.
We will interested in the transformation low for
the first five coefficients:

\proclaim Lemma 7.1. The action 
$\sigma\circ f_{\nu,\rho}\circ\sigma^{-1}$
of the group ${\rm Diff}({\bf RP}^1)$
on the quotient-module
${\cal D}^k_{\nu,\rho}/{\cal D}^{k-5}_{\nu,\rho}$
is given by:
$$
\matrix{
f(\bar a_k)  &=&  f^*(\bar a_k)\hfill\cr
\noalign{\smallskip}
f(\bar a_{k-1})  &=& f^*(\bar a_{k-1})\hfill\cr 
\noalign{\smallskip}
f(\bar a_{k-2})  &=&
f^*(\bar a_{k-2})+
\beta_{k-2}^kS(f^{-1})\bar a_k\hfill\cr
\noalign{\bigskip}
f(\bar a_{k-3} ) &=& f^*(\bar a_{k-3})+
\beta_{k-3}^{k-1}S(f^{-1})\bar a_{k-1}+
\beta_{k-3}^k{\cal T}_{\rho-\nu-k}(f^{-1})(\bar a_k)\hfill\cr
\noalign{\bigskip}
f(\bar a_{k-4}) &=& f^*(\bar a_{k-4})+
\beta_{k-4}^{k-2}S(f^{-1})\bar a_2+
\beta_{k-4}^{k-1}{\cal T}_{\rho-\nu-k+1}(f^{-1})(\bar a_{k-1})
\hfill\cr
\noalign{\smallskip}
&&\;\;\;\;\;\;\;\;\;\;\;\;\;\;\;\;
+\beta_{k-4}^k{\cal U}_{\rho-\nu-k}(f^{-1})(\bar a_k)\hfill\cr
}
\eqno{(\ref{act}'')}
$$
where 
$f^*(\bar a_i)=
\bar a_i\circ f^{-1}\cdot({f^{-1}}^{\prime})^{\nu-\rho-i}$,
and $\beta_i^j$ are some constants
depending on $\nu$ and $\rho$.
\par

\vskip 0,3cm
\noindent
{\bf 7.3}.
Recall that the cohomology group
$H^1(G;{\rm Hom}(A,B))$,
where $G$ is a group, $A$ and $B$ are $G$-modules,
classifies nontrivial extensions of
$G$-modules:
$$
0\to A\to{\cal E}\to B\to0.
$$
Given a 1-cocycle 
${\cal C}:G\to
{\rm Hom}(A,B)$,
one has a 
$G$-module structure on
${\cal E}=A\oplus B$ 
defined by: 
$$
\rho_g(\phi ,\psi):=
(g(\phi) ,g(\psi) +\gamma {\cal C}(g)(g(\phi)) ),
$$
where $\gamma$ is an arbitrary constant. 
Moreover, the condition 
$\rho_f\circ \rho_g=\rho_{f\circ g}$ is equivalent to the fact
that ${\cal C}$ is a 1-cocycle.

\vskip 0,3cm

Let us realize the extensions
$
{\cal E}_{\lambda,\lambda+2}={\cal F}_{\lambda}\oplus {\cal F}_{\lambda+2}
$
and
$
{\cal E}_{\lambda,\lambda+3}={\cal F}_{\lambda}\oplus {\cal F}_{\lambda+3}
$
defined by 
the cocycles ${\cal S}_{\lambda}$ and ${\cal T}_{\lambda}$
respectively,
as some modules of differential operators.

\vskip 0,3cm
\noindent
{\bf 7.4 Examples}. 1) If $\lambda\not=0,-\frac{1}{2},-1$, then
the submodule of 
${\cal D}^2_{\nu,\nu+\lambda+2}$ consisting of differential operators:
$$
A=a_2(x)\frac{d^2}{dx^2}-
\frac{2\nu+1}{\lambda}a_2'(x)\frac{d}{dx}+a_0(x)
$$
is isomorphic to 
the module ${\cal E}_{\lambda,\lambda+2}$;

\vskip 0,3cm

2) If $\lambda\not=0,-\frac{1}{2},-1,-\frac{3}{2},-2$, then
the submodule of 
${\cal D}^3_{\nu,\nu+\lambda+3}$,
where $\nu$ is determined by the
condition:
$$
3\nu^2+3\nu(\lambda+2)+\lambda+2=0
$$
and the differential operators are of the form:
$$
A=a_3(x)\frac{d^3}{dx^3}-
3\frac{\nu+1}{\lambda}a_3'(x)\frac{d^2}{dx^2}+
3\frac{(\nu+1)(2\nu+1)}{\lambda(2\lambda+1)}a_3''(x)\frac{d}{dx}+
a_0(x)
$$
is isomorphic to 
the module ${\cal E}_{\lambda,\lambda+3}$.

\section{Appendix}

It follows from Theorem 1 that there exists two
1-cocycles on $\Diff({\bf RP}^1)$
vanishing on $PSL(2{\bf R})$, with
values on ${\cal D}_{\lambda,\lambda+5}$
for the particular values of $\lambda=0$ and $-4$.

Let us give the explicit formul{\ae}
for non-trivial cocycles on
$\Diff({\bf RP}^1)$ with values in
${\cal D}^3_{0,5}$ and ${\cal D}^3_{-4,1}$:
$$
V_0(f)=
S(f)\frac{d^3}{dx^3}
-\frac{3}{2}S(f)^{\prime}\frac{d^2}{dx^2}
+\left(\frac{3}{10}S(f)^{\prime\prime}-
\frac{4}{5}S(f)^2\right)\frac{d}{dx}
$$
and
$$
\matrix{
\displaystyle V_{-4}(f)=
S(f)\frac{d^3}{dx^3}
+\frac{9}{2}S(f)^{\prime}\frac{d^2}{dx^2}
+\left(\frac{63}{10}S(f)^{\prime\prime}-
\frac{4}{5}S(f)^2\right)\frac{d}{dx}\;\;
\hfill\cr\noalign{\smallskip}
\;\;\;\;\;\;\;\;\;\;\;\;\;\;\;\;\;\;\;\;\;\;\;\;\;\;\;\;\;\;
\;\;\;\;\;\;\;\;\;\;\;\;\;\;\;\;\;
\displaystyle +\frac{14}{5}S(f)^{\prime\prime\prime}-
\frac{8}{5}S(f)^{\prime}S(f)
}
$$
respectively.

The proof is straightforward.

\vskip 0,5cm

{\it Acknowledgments}. It is a pleasure to acknowledge 
numerous fruitful discussions 
with C. Duval, B.Kostant, P. Lecomte and E. Mourre.

\vfill\eject


\vskip 1cm

\noindent
Sofiane BOUARROUDJ\\
{\small CMI, Universit\'e de Provence}\\
{\small 39 rue Juliot Curie}\\
{\small 13453 Marseille, France}\\
\\
Valentin OVSIENKO\\
{\small C.N.R.S., C.P.T.}\\
{\small  Luminy-Case 907}\\
{\small  F-13288 Marseille Cedex 9, France}

\end{document}